# On random walks in random scenery

## F. M. Dekking[1] and P. Liardet[2]

*Delft University of Technology and Université de Provence*

**Abstract:** This paper considers 1-dimensional generalized random walks in random scenery. That is, the steps of the walk are generated by an arbitrary stationary process, and also the scenery is a priori arbitrary stationary. Under an ergodicity condition—which is satisfied in the classical case—a simple proof of the distinguishability of periodic sceneries is given.

## 1. Introduction

Random walks in random scenery have been studied by Mike Keane for quite some time (see [2] for his most recent work). In fact, he and Frank den Hollander were pioneers in this exciting area. Around 1985 they formulated a conjecture about "recovery of the scene" by a simple random walker. A weaker form of this: "distinguishability of two scenes", was proven by Benjamini and Kesten ([1]). Since then there has been a lot of action in this field, especially by Matzinger and co-workers. We just mention the recent paper [5].

In the following we will introduce generalized random walks in random scenery, and analyse them from a dynamical point of view. This gives us in Section 2 a general scenery recovery result on the level of measures, from which we deduce in a simple way in Section 3 a proof for the distinguishability of periodic sceneries.

We consider a random walker on the integers. The integers are coloured by colours from an alphabet $C$. This is the scenery. At time $n$ the walker records the scenery at his position, this yields $r_n$ from $C$. To formalize somewhat more, let the random walk be described by a measure $\mu$ on the Borel sets of $\Omega$, where

$$\Omega = \{\omega = (\omega_n)_{n \in \mathbb{Z}} : \omega_n \in J \text{ for all } n\}.$$

Here the set $J$ of the possible steps of the walk will simply be $\{-1, +1\}$, or somewhat more general $\{-1, 0, +1\}$. Although often a single scenery $x = (x_k)_{k \in \mathbb{Z}}$ is considered, it is useful to consider $x$ as an element of the shift space $X = C^{\mathbb{Z}}$ with shift map $T : X \to X$, equipped with some ergodic $T$-invariant measure $\lambda$, which we will call the *scenery measure*. We then consider $x$ picked according to the measure $\lambda$. The colour record $\varphi_x$ of $x$ can be written as a map $\varphi_x : \Omega \to X$:

$$\varphi_x(\omega) = (r_n(\omega, x))_{n \in \mathbb{Z}},$$

where in line with the description above, one has for $n \geq 1$

$$r_n(\omega, x) = (T^{\omega_0 + \cdots + \omega_{n-1}} x)_0.$$

---

[1]Thomas Stieljes Institute for Mathematics and Delft University of Technology Faculty EEMCS, Mekelweg 4, 2628 CD Delft, The Nethrlands, e-mail: F.M.Dekking@math.tudelft.nl

[2]Université de Provence, Centre de Mathématiques et Informatique (CMI), 39 rue Joliot-Curie, F-13453 Marseille cedex 13, France, e-mail: liardet@cmi.univ-mrs.fr

This definition is completed by putting $r_0(\omega, x) = x_0$, and for $n < 0$:

$$r_n(\omega, x) = (T^{-\omega_{-1} - \omega_{-2} \cdots - \omega_n} x)_0.$$

The dynamics of the whole process is well described by a skew product transformation $T_{\Omega \times X}$ on the product space $\Omega \times X$ defined by

$$T_{\Omega \times X}(\omega, x) = (\sigma\omega, T^{\omega_0} x),$$

where $\sigma$ denotes the shift map on $\Omega$.

Let us now look at the colour records of *all* $x$; we define the *global recording map* $\Phi : \Omega \times X \to X$ by

$$\Phi(\omega, x) = (r_n(\omega, x))_{n \in \mathbb{Z}}.$$

**Lemma 1.** *The map $\Phi$ is equivariant, that is, $\Phi \circ T_{\Omega \times X} = T \circ \Phi$.*

*Proof.* One way:

$$\begin{aligned}
\Phi \circ T_{\Omega \times X}(\omega, x) &= \Phi((\sigma\omega, T^{\omega_0} x)) = (r_n(\sigma\omega, T^{\omega_0} x))_{n \in \mathbb{Z}} \\
&= ((T^{\omega_1 + \cdots + \omega_n} T^{\omega_0} x)_0) = (r_{n+1}(\omega, x))_{n \in \mathbb{Z}}.
\end{aligned}$$

The other way:

$$T \circ \Phi(\omega, x) = T((r_n(\omega, x))_{n \in \mathbb{Z}}) = (r_{n+1}(\omega, x))_{n \in \mathbb{Z}}. \qquad \square$$

Clearly product measure $\mu \times \lambda$ is preserved by $T_{\Omega \times X}$. We will be particularly interested in the image of $\mu \times \lambda$ under the global recording map $\Phi$, which we denote $\rho$:

$$\rho = (\mu \times \lambda) \circ \Phi^{-1}.$$

We call $\rho$ the *global record measure*. It follows from Lemma 1 that $\rho$ is invariant for $T$. Moreover, $\rho$ will be ergodic when $T_{\Omega \times X}$ is ergodic for $\mu \times \lambda$. In the classical case were $\mu$ is product measure this is guaranteed by Kakutani's random ergodic theorem. In this case, when $\lambda$ and $\lambda'$ are two scenery measures, and $\rho = (\mu \times \lambda) \circ \Phi^{-1}$ and $\rho' = (\mu \times \lambda') \circ \Phi^{-1}$ are the corresponding global record measures, then either $\rho = \rho'$ or $\rho \perp \rho'$.

The colour record $\varphi_x$ of a scenery $x$ induces the *record measure* $\rho_x$ defined by

$$\rho_x = \mu \circ \varphi_x^{-1}.$$

Following [4] we call the two sceneries $x$ and $y$ *distinguishable* if $\rho_x \perp \rho_y$. The following lemma shows that global distinguishability carries over to local distinguishability.

**Lemma 2.** *Let $\lambda$ and $\lambda'$ be two scenery measures with corresponding global record measures $\rho$ and $\rho'$. Then $\rho \perp \rho'$ implies that $\rho_x \perp \rho_y$ for $\lambda \times \lambda'$ almost all $(x,y)$.*

*Proof.* By Fubini's theorem, and recalling that $\Phi(\omega, x) = \varphi_x(\omega)$,

$$\begin{aligned}
\rho(E) &= \int_X \int_\Omega 1_E \circ \Phi(\omega, x) \, \mathrm{d}\mu(\omega) \, \mathrm{d}\lambda(x) \\
&= \int_X \mu(\varphi_x^{-1} E) \, \mathrm{d}\lambda(x).
\end{aligned}$$

So $\rho = \int \rho_x \, \mathrm{d}\lambda(x)$. Hence if $E$ is a Borel set with the property that $\rho(E) = 1$ and $\rho'(E^c) = 1$, then $\rho_x(E) = 1$ for $\lambda$-almost $x$, and $\rho_y(E^c) = 1$ for $\lambda'$-almost $y$. $\qquad \square$

## 2. Reconstructing the scenery measure

Here we consider the case of a generalized random walk with steps $J = \{-1, 0, +1\}$ given by an ergodic stationary measure $\mu$ on $\Omega = J^{\mathbb{Z}}$. For ease of notation we rename $J$ to $\{\mathtt{L}, \mathtt{H}, \mathtt{R}\}$. To simplify the exposition we will assume that there is no holding ($\mu([\mathtt{H}]) = 0$), and show at the appropriate moment that this restriction can trivially be removed.

There is a basic difference between symmetric walks and asymmetric walks in the reconstruction of the scenery. We call $\mu$ *symmetric* if for each word $w$

$$\mu[w] = \mu[\overline{w}].$$

Here $\overline{w}$ denote the mirror image of $w$, that is, the word obtained from $w$ by replacing $\mathtt{R}$ by $\mathtt{L}$ and $\mathtt{L}$ by $\mathtt{R}$. Since he does not know left from right, a symmetric walker can only reconstruct a scenery $x$ up to a reflection. This will result in two theorems, one for the asymmetric, and one for the symmetric case.

Let $\lambda$ be a scenery measure. We give a few examples of calculation of the $\rho$-probabilities of cylinder sets. Let $\mathcal{W}_n$ be the set of all words $w = w_1 \ldots w_n$ of length $n$ over the (colour) alphabet $\{0, 1\}$. For $w \in \mathcal{W}_n$ we let $[w]$ denote the cylinder

$$[w] = \{x \in X : x_0 \ldots x_{n-1} = w\},$$

and we will abbreviate $\rho([w])$ to $\rho[w]$. We will use the same type of notations and conventions for $\lambda$ and $\mu$. It is clear, using the stationarity of $\lambda$, that for instance

$$\rho[001] = \mu[\mathtt{RR}]\lambda[001] + \mu[\mathtt{LL}]\lambda[100],$$

and slightly more involved

$$\rho[000] = (\mu[\mathtt{RL}] + \mu[\mathtt{LR}])\lambda[00] + (\mu[\mathtt{RR}] + \mu[\mathtt{LL}])\lambda[000],$$

and

$$\rho[0001] = \mu[\mathtt{RLL}]\lambda[100] + \mu[\mathtt{LRR}]\lambda[001] + \mu[\mathtt{RRR}]\lambda[0001] + \mu[\mathtt{LLL}]\lambda[1000].$$

In the sequel we shall denote the word $\mathtt{R} \ldots \mathtt{R}$, $N$ times repeated, as $\mathtt{R}^N$. Note how with each appearance of a word $w$ on the right side also the reversed word $\overleftarrow{w}$ appears, defined by $\overleftarrow{w} = w_n \ldots w_1$ if $w = w_1 \ldots w_n$. Words $w$ that satisfy $w = \overleftarrow{w}$ are called *palindromes*. Now let us put all the words $w$ from $\cup_{1 \le k \le n} \mathcal{W}_k$ in some fixed order, taking care that their lengths are non-decreasing and that for a fixed $k$ we first take all palindromes, and then all non-palindromes in pairs $(w, \overleftarrow{w})$. Let $V_n(\rho)$ and $V_n(\lambda)$ denote the vectors of length $(2^{n+1} - 2)$ containing the real numbers $\rho[w]$ respectively $\lambda[w]$ in the chosen order. For example,

$$V_2^{\mathrm{T}}(\rho) = (\rho[0], \rho[1], \rho[00], \rho[11], \rho[01], \rho[10]).$$

In general, if $w$ is a word of length $N + 1$, then $\rho[w]$ is obtained as a sum of products $\mu[u]\lambda[v]$, where the length of $v$ is *at most* $N + 1$, and length $N + 1$ only occurs when the walker makes no turns, i.e., when $u = \mathtt{R}^N$ or $u = \mathtt{L}^N$. Moreover, if $w$ is a palindrome, then there is one maximal length term $(\mu[\mathtt{R}^N] + \mu[\mathtt{L}^N])\lambda[w]$, and if $w$ is not a palindrome, then there are two maximal length terms $\mu[\mathtt{R}^N]\lambda[w]$, and $\mu[\mathtt{L}^N]\lambda[\overleftarrow{w}]$. This observation shows that there exists an almost lower triangular $(2^{n+1} - 2) \times (2^{n+1} - 2)$ matrix $A_n$ such that

$$V_n(\rho) = A_n V_n(\lambda).$$

Here 'almost lower triangular' means that $A_n$ has the form

$$
\begin{pmatrix}
\square & 0 & 0 & 0 & 0 & 0 \\
* & \square & 0 & 0 & 0 & 0 \\
* & * & \square & 0 & 0 & 0 \\
* & * & * & \ldots & 0 & 0 \\
* & * & * & * & \square & 0 \\
* & * & * & * & * & \square
\end{pmatrix},
$$

where (at palindrome entries) $\square$ is a $1 \times 1$ matrix $\mu[\mathrm{R}^N] + \mu[\mathrm{L}^N]$, and (at non-palindrome pairs) $\square$ is a $2 \times 2$ matrix of the form

$$
\begin{pmatrix}
\mu[\mathrm{R}^N] & \mu[\mathrm{L}^N] \\
\mu[\mathrm{L}^N] & \mu[\mathrm{R}^N]
\end{pmatrix}.
$$

With simple linear algebra we find that $A_n$ is non-singular if and only if

$$
\mu[\mathrm{R}] \neq \mu[\mathrm{L}], \ldots, \mu[\mathrm{R}^N] \neq \mu[\mathrm{L}^N], \ldots
$$

Let us call a generalized random walk given by $\mu$ *strongly asymmetric* if all these inequalities hold. For instance, if $\mu$ is a stationary Markov chain given by a $2 \times 2$ transition matrix $(p_{s,s})$, then $\mu$ is strongly asymmetric if and only if $p_{\mathrm{R,R}} \neq p_{\mathrm{L,L}}$.

Note that when $\mu[\mathrm{H}] > 0$, then only some *sub*-diagonal elements of $A_n$ will change from 0 to a positive value. We therefore obtained the following result.

**Theorem 1.** *For strongly asymmetric generalized random walk with holding the scenery measure $\lambda$ can be reconstructed from $\rho$.*

What remains is the symmetric walker case. Then in general $\lambda$ can not be reconstructed from $\rho$. However, often we can reconstruct the *reversal symmetrized measure* $\check{\lambda}$ defined for each word $w$ by

$$
\check{\lambda}[w] = \frac{1}{2}\Big(\lambda[w] + \lambda[\overleftarrow{w}]\Big).
$$

For symmetric $\mu$ the equation for, e.g., $\rho[0001]$ becomes

$$
\rho[0001] = 2\mu[\mathrm{LRR}]\check{\lambda}[001] + 2\mu[\mathrm{RRR}]\check{\lambda}[0001].
$$

In general, if $w$ is a word of length $N+1$, then $\rho[w]$ is obtained as a sum of products $\mu[u]\check{\lambda}[v]$, where the length of $v$ is *at most* $N+1$, and length $N+1$ only occurs when $u = \mathrm{R}^N$ or $u = \mathrm{L}^N$. Moreover, now there is for all words $w$ one term $2\mu[\mathrm{R}^N]\check{\lambda}[w]$ for the $v = w$ with maximal length. So this time we obtain the existence of a $(2^{n+1} - 2) \times (2^{n+1} - 2)$ lower triangular matrix $A_n$ such that

$$
V_n(\rho) = A_n V_n(\check{\lambda}).
$$

Let us call $\mu$ *straightforward* if arbitrary long words of $\mathrm{R}$'s have positive probability to appear. Then the diagonal elements of $A_n$ are positive for each $n$, and we obtain the following.

**Theorem 2.** *For straightforward symmetric generalized random walk with holding the reversal symmetrized scenery measure $\check{\lambda}$ can be reconstructed from $\rho$.*

## 3. Distinguishing periodic sceneries

In this section we shall derive more general results with more simple proofs than in [3]. It is shown there that for asymmetric simple random walk with holding any two periodic sceneries $x$ and $y$ which are not translates of each other can be distinguished by their scenery records, i.e., $\rho_x \perp \rho_y$. Our result is

**Theorem 3.** *Any strongly asymmetric generalized random walk with holding can distinguish two periodic sceneries that are not translates of each other, provided that their global record measures are ergodic.*

*Proof.* Let us write $x \overset{\mathrm{t}}{\sim} y$ if $x$ and $y$ are translates of each other, i.e., for some $k$ one has $y = T^k x$. Let $\mathrm{Per}(x)$ be the period of $x$, i.e. $p = \mathrm{Per}(x)$ is the smallest natural number such that $T^p x = x$. Let $\lambda$ be the scenery measure generated by $x$, i.e., denoting point measure in $z$ by $\delta_z$,

$$\lambda = \frac{1}{\mathrm{Per}(x)} \sum_{k=0}^{\mathrm{Per}(x)-1} \delta_{T^k x}.$$

The scenery measure generated by $y$ is denoted as $\lambda'$. Now suppose that $\rho_x$ is *not* orthogonal to $\rho_y$. Then, since $\lambda$ and $\lambda'$ are discrete, it follows from Lemma 2 that also $\rho$ and $\rho'$ are not orthogonal. But since these measures are ergodic, they must be equal. From Theorem 1 it then follows that also $\lambda = \lambda'$. This implies that $x \overset{\mathrm{t}}{\sim} y$, by the discreteness of $\lambda$ and $\lambda'$. Indeed, equality of these measures yields that $\delta_{T^k x} = \delta_{T^j y}$ for some $k$ and $j$, and hence that $x \overset{\mathrm{t}}{\sim} y$. $\qquad \square$

For a symmetric (generalized) random walk it is impossible to distinguish a sequence $x$ from its reflection $\overset{\leftarrow}{x}$, defined by $\overset{\leftarrow}{x}_k = x_{-k}$. So let us call $x$ and $y$ equivalent, and we denote $x \sim y$, if $y$ can be obtained from $x$ by translation and/or reflection.

**Theorem 4.** *Any straightforward symmetric generalized random walk with holding can distinguish two periodic sceneries that are not equivalent, provided that their global record measures are ergodic.*

*Proof.* The proof follows the same path as the proof of Theorem 3, using Theorem 2 instead of Theorem 1. The only other difference now is that the measure $\tilde{\lambda}$ is a mixture of point measures in $T^k x$ and in $T^j \overset{\leftarrow}{x}$. But then equality of $\tilde{\lambda}$ and $\tilde{\lambda}'$ implies that $y$ must be a translate, or the reflection of a translate of $x$, i.e., $x \sim y$. $\qquad \square$

### Acknowledgment

We are grateful, in chronological order, to Jeff Steif and Frank den Hollander for their useful comments on earlier versions of this paper.